\documentclass[12pt,reqno]{amsart}
\usepackage{amssymb,latexsym,amsmath,mathrsfs, xypic} 
\usepackage[pdftex, colorlinks=false, pdfborder={1 1 1}, pdfborderstyle={/S/U/W 1}]{hyperref} 

\newcommand{\C}{\mathbb{C}} 
\newcommand{\N}{\mathbb{N}}

\newcommand{\Z}{\mathbb{Z}}

\newcommand{\ubf}{\mathbf{u}}
\newcommand{\zbf}{\mathbf{z}}

\newcommand{\Od}{{\mathcal{O}_2}}

\newcommand{\Td}{{\mathcal{T}_2}}

\newcommand{\ot}{\otimes}
\newcommand{\pf}{\noindent {\mbox{\textit{Proof}. }} }
\newcommand{\ie}{\textit{i.e.}\,\,} 
\newcommand{\eg}{\textit{e.g.\/}\ }
\newcommand{\cst}{C$^*$}

\newtheorem{thm}{Theorem}
\newtheorem{cor}[thm]{Corollary}
\newtheorem{lem}[thm]{Lemma}
\newtheorem{prop}[thm]{Proposition}

\newtheorem{rems}[thm]{Remarks}

\numberwithin{equation}{section}
\title{A properly infinite \cst-algebra \\which is not K$_1$-injective}

\author{Etienne Blanchard}



\begin{document} 
\begin{abstract} 
We construct in this note a unital properly infinite \cst-algebra which is not K$_1$-injective. 
\end{abstract} 

\maketitle
\vspace{-15pt} ${}$\hfill In memory of Uffe Haagerup

\section{Notations} 
\indent The \cst-algebra $\Od$ defined by J. Cuntz in \cite{Cun77} 
is the unital \cst-algebra generated by two isometries $s_1, s_2$ 
satisfying the equality $s_1^{}s_1^*+s_2^{}s_2^*=1$. 
\\ \indent 
We call $\Td$ the extension of $\Od$ by the \cst-algebra of compact operators acting on 
an infinite dimensional Hilbert space. 
This universal unital \cst-algebra is generated by two isometries $v_1, v_2$ satisfying 
the inequality $v_1^{}v_1^*+v_2^{}v_2^*\leq 1$ (\cite{Cun77}). 

\medskip 
A nonzero unital \cst-algebra $A$ is said to be \textit{properly infinite} if and only if 
there exists an homomorphism of unital \cst-algebra $\Td\to A$ (see \cite[Proposition 2.1]{Ror04}). 
\\ \indent 
A nonzero unital \cst-algebra $A$ is said to be \textit{K$_1$-injective} if and only if 
any unitary $u\in {A}$ satisfies $[u]=[1_A]$ in $K_1(A)$ only if it
belongs to the connected component $\mathcal{U}^0(A)$ of the unit $1_A$ 
in the group $\mathcal{U}(A)$ of unitaries in $A$ (see \eg \cite[definition 2.6]{BRR08}). 
\\ \indent 
We write $u_1\sim_h u_2$ when two unitaries $u_1, u_2$ in $\mathcal{U}(A)$ satisfy 
$u_1^{}\cdot u_2^*\in\mathcal{U}^0(A)$. 

We answer in Corollary~\ref{deformpropinf} several questions 
from \cite{BRR08}, \cite{Blan09}, \cite{Blan16} 
on the stability of proper infiniteness under continuous deformations. 

\section{Main result} 
Let $\mathbf{u}$ be the canonical unitary generator of the \cst-algebra $C^*(\Z)$ and 
denote by $\imath_0$ and $\imath_1$ the two canonical $\ast$-embeddings 
of the \cst-algebra $\Od$ in the full unital free product $\Od\ast_\C\Od$. 
Then one has the following description 
for the K-theory of the unital free product $\Od\ast_\C\Od\,$. 

\begin{lem}\label{lem2.1} 
1) $K_0(\Od\ast_\C\Od)=0$ and $K_1(\Od\ast_\C\Od)=\Z$. 
\\ 2) There is an isomorphism of unital \cst-algebra $\Od\ast_\C C^*(\Z)\cong\Od\ast_\C\Od$. 
\\ 3) The unitary $\mathbf{u}$ 
has a generating image in $K_1(\Od\ast_\C C^*(\Z))$. 
\end{lem} 
\pf 
1) As noticed by B. Blackadar in \cite{Blac07}, 
Corollary~2.6 of \cite{Ger97} implies 
the exactness of the six-term cyclic sequence of topological K-theory group 
$$\begin{array}{cccccc}
K_0(\C)=\Z&\mathop{\longrightarrow}&K_0(\Od\oplus \Od)=0\oplus0&\mathop{\longrightarrow}&
K_0(\Od\ast_\C\Od)\\
\uparrow&&&&\downarrow\\
K_1( \Od\ast_\C\Od)&\longleftarrow &K_1(\Od\oplus\Od)=0\oplus0&
\longleftarrow&K_1(\C)=0
\end{array}$$ 
2) The unital free product $\Od\ast_\C C^*(\Z)$ is generated by 
the two isometries $s_1, s_2$ and the unitary $\mathbf{u}$ 
whereas the unital free product $\Od\ast_\C\Od$ is generated by 
the four isometries $\imath_0(s_1), \imath_0(s_2),\imath_1(s_1), \imath_1(s_2)$. 
The unique \hbox{\cst-morphism} $\sigma: \Od\ast_\C C^*(\Z)\to\Od\ast_\C\Od$ with 
\begin{center} 
{$\sigma(s_k)=\imath_0(s_k)\;$ ($k=1, 2$)\quad and\quad 
$\sigma(\mathbf{u})=\imath_1(s_1^{})\imath_0(s_1^*)+\imath_1(s_2^{})\imath_0(s_2^*)$} 
\end{center} 
is an isomorphism satisfying $\sigma^{-1}\big(\imath_0(s_k)\bigr)=s_k$ and 
$\sigma^{-1}\big(\imath_1(s_k)\bigr)=\mathbf{u}\cdot s_k$ for $k= 1, 2$. 

\noindent 3) The unitary $\mathbf{u}$ generates a canonical copy of the \cst-algebra $C^*(\Z)$ 
inside the unital free product $\Od\ast_\C C^*(\Z)$ and 
corollary~2.6 of \cite{Ger97} entails the exactness of the cyclic sequence 
$$\begin{array}{ccccc}
K_0(\C)=\Z&\mathop{\longrightarrow}&K_0(\Od\oplus C^*(\Z))=0\oplus\Z&\mathop{\longrightarrow}&
K_0(\Od\ast_\C C^*(\Z))=0\\
\uparrow&&&&\downarrow\\
K_1( \Od\ast_\C C^*(\Z))=\Z&\longleftarrow &K_1(\Od\oplus C^*(\Z))=0\oplus\Z&
\longleftarrow&K_1(\C)=0
\end{array}$$
The injectivity of the upper left arrow implies that 
the lower left arrow is an isomorphism and 
the class $[\mathbf{u}]$ is a generator of the group K$_1(\Od\ast_\C C^*(\Z))\cong\Z$. \qed

\medskip\indent 
Denote by $\jmath_0$ and $\jmath_1$ the two canonical embeddings 
of the \cst-algebra $\Td$ in the unital free product $\Td\ast_\C\Td$. 
As the projection $1-v_1^{}v_1 ^*$ is properly infinite and full in $\Td$, 
the two projections $\jmath_0(1-v_1^{}v_1 ^*), \jmath_1(1-v_1^{}v_1 ^*)$ 
are properly infinite and full in $\Td\ast_\C\Td$ and 
there exists (Lemma~2.4 in \cite{BRR08}) a unitary $\tilde{u}\in\mathcal{U}(\Td\ast_\C\Td)$ with 
$\jmath_1(v_1)=\tilde{u}\cdot \jmath_0(v_1)$. 
 \\ \indent 
Define also the K$_1$-trivial unitary $\mathring{u}$ 
in the quotient $\Od\ast_\C\Od=\Od\ast_\C C^*(\Z)$ by 
$$\begin{array}{cl}\mathring{u}\hspace{-0,2cm}&=
[\imath_1(s_1)\imath_0(s_1^*)+\imath_1(s_2)\imath_0(s_2^*)]\cdot
[\imath_0(s_1^{}s_1^*)+
\imath_0(s_2)[\imath_0(s_1)\imath_1(s_1^*)+\imath_0(s_2)\imath_1(s_2^*)]\imath_0(s_ 2^*)] \\
&=\mathbf{u}\cdot(s_1^{}s_1^*+s_2^{}\mathbf{u}s_2^*)^{-1}\;. 
\end{array}$$ 
\indent Then the unital free product $\Td\ast_\C\Td$ enjoys the following properties. 

\begin{prop}\label{prop2} 
1) $K_0(\Td\ast_\C\Td)=\Z$ and $K_1(\Td\ast_\C\Td)=0$. \\
2) The direct sum $\tilde{u}\oplus 1_{\Td\ast_\C\Td}$ belongs to the connected component 
$\mathcal{U}^0(M_2(\Td\ast_\C\Td))\,$. \\ 
3) The unitary $\tilde{u}$ is in $\mathcal{U}^0(\Td\ast_\C\Td)$ if and only if 
the unitary $\mathring{u}$ is in $\mathcal{U}^0(\Od\ast_\C C^*(\Z))$. \\ 
4) Any unital $\ast$-representation $\Theta$ of the unital free product $\Od\ast_\C C^*(\Z)$ 
which factors through the tensor product $\Od\ot C^*(\Z)$ satisfies 
$\Theta(\mathring{u})\in\mathcal{U}^0\left(\Theta(\Od\ast_\C C^*(\Z))\right)$. \\
5) The unitary $\mathring{u}$ is not in the connected component 
$\mathcal{U}^0\left(\Od\ast_\C C^*(\Z)\right)$. \\ 
6) The properly infinite \cst-algebra $\Td\ast_\C\Td$ is not $K_1$-injective. 
\end{prop}
\pf 1) The suspension $C_0((0,1))\otimes\Td\ast_\C\Td$ is $KK$-equivalent to 
the cone $D=\{(f_0, f_1)\in C_0( (0, 1], \Td\oplus\Td), f_0(1)=f_1(1)\in\C\cdot 1_\Td\}$ 
(Theorem 2.2 in \cite{Ger97}) and that cone is $KK$-equivalent 
to the \cst-algebra $\{(f_0, f_1)\in C((0, 1])^2, f_0(1)=f_1(1)\}\cong C_0((0, 2))$. 
Thus, $K_i(\Td\ast_\C\Td)=K_{1-i}(D)=K_{1-i}(C_0(0, 2))$ for $i=0, 1$. 

\medskip\noindent
2) The two projections $0\oplus1$ and $1\oplus0$ are properly infinite and full 
in $M_2(\Td\ast_\C\Td)$. 
Hence, the unitary $\tilde{u}\oplus 1$ belongs to $\mathcal{U}^0(M_2(\Td\ast_\C\Td))$ 
(see \cite[exercise~8.11]{LLR00}). 

\medskip\noindent
3) Call $\pi$ the unique quotient morphism of unital \cst-algebra 
from the full unital free product $\Td\ast_\C\Td$ to $\Od\ast_\C\Od$ 
such that $\pi(\jmath_i(v_k))=\imath_i(s_k)$ for all $i=0, 1$ and $k=1, 2$. 
\indent Then $\pi(\tilde{u})s_1=\pi(\tilde{u}\jmath_0(v_1))=\pi(\jmath_1(v_1))=\mathbf{u}s_1$ and 
the product $\pi(\tilde{u})\cdot(\mathring{u})^*$ satisfies 
$$\begin{array}{rll}
\pi(\tilde{u})\cdot(\mathring{u})^*\cdot\ubf s_1^{}s_1^*\ubf^*
&=\pi(\tilde{u})\cdot s_1^{}\cdot s_ 1^*\ubf^*&\\
&=\ubf s_1^{}\cdot s_1^*\ubf^*&\quad\mathrm{in}\quad \Od\ast_\C C^*(\Z)\,.
\end{array}$$ 
\indent Accordingly, the product $\pi(\tilde{u})\cdot(\mathring{u})^*$ is 
a K$_1$-trivial unitary in $\Od\ast_\C C^ *(\Z)$ 
which commutes to the projection $\ubf s_1^{}s_1^*\ubf^*$. 
As both $\ubf s_1^{}s_1^*\ubf^*$ and $1- \ubf s_1^{}s_1^*\ubf^*$ 
are properly infinite full projections in $\Od\ast_\C C^ *(\Z)$, 
Lemma~2.4 of \cite{BRR08} implies that 
$\pi(\tilde{u})\cdot(\mathring{u})^*$ belongs to $\mathcal{U}^0(\Od\ast_\C C^ *(\Z))$ 
and so $\pi(\tilde{u})=\pi(\tilde{u})\cdot(\mathring{u})^*\cdot\mathring{u}\sim_h\mathring{u}$ 
in $\mathcal{U}(\Od\ast_\C C^ *(\Z))$. 

\noindent -- If $\tilde{u}\sim_h 1_{\Td\ast_\C\Td}$ in $\mathcal{U}(\Td\ast_\C\Td)$, then 
$\mathring{u}\sim_h\pi(\tilde{u})\sim_h \pi(1_{\Td\ast_\C\Td})=
1_{\Od\ast_\C\Od}$ in $\mathcal{U}(\Od\ast_\C\Od) $. 

\noindent -- If $\mathring{u}\sim_h 1_{\Od\ast_\C\Od}$ in $\mathcal{U}(\Od\ast_\C\Od)$, 
then the homotopic unitary $\pi(\tilde{u})$ also belongs to 
the connected component $\mathcal{U}^0(\Od\ast_\C\Od)$ and 
its lift $\tilde{u}$ in $\mathcal{U}(\Td\ast_\C\Td)$ can be continuously connected 
to the unit $1_{\Td\ast_\C\Td}$ (see \eg Lemma~2.1.7 in \cite{LLR00}). 

\medskip\noindent4) Let $\sigma: \Od\ast_\C C^*(\Z)=C^*\langle s_1, s_2, \ubf\rangle\to\Od\otimes C^*(\Z)=
C(S^1, \Od)$ be the unique \cst-epimorphism satisfying the relations 
\begin{center} 
$\sigma(s_k)=s_k\otimes 1\;$ ($k=1, 2$)\quad and\quad 
$\sigma(\ubf)=1\otimes\ubf$\,. 
\end{center}
The two isometries $V_1= (s_1^{}s_1^*+s_2^{}s_1^{}s_2^*)\otimes1$ and 
$V_2=s_2^{}s_2^{}\otimes1$ generate a unital copy of $\Od$ in $\Od\otimes C^*(\Z)$
and we have in $\mathcal{U}(\Od\otimes C^*(\Z))$ the sequence of homotopies 
$$\begin{array}{rll}
\sigma(\mathring{u})\hspace{-0,3cm}&=\sigma(\ubf(s_1^{}s_1^*+s_2^{}\ubf^*s_2^*))\\
&=s_1^{}s_1^*\otimes\ubf+
(1\otimes\ubf)(s_2^{}\otimes 1)(1\otimes\ubf^*)(s_2^*\otimes 1)\\
&=s_1^{}s_1^*\otimes\ubf+s_2^{}s_2^*\otimes1\\
&=s_1^{}s_1^*\otimes\ubf+s_2^{}s_1^{}s_1^*s_2^*\otimes 1+
s_2^{}s_2^{}s_2^*s_2^*\otimes 1&\mathrm{because}\,1_\Od=s_1^{}s_1^*+s_2^{}s_2^*\\
&=V_1^{}\sigma(\mathring{u})V_1^*+V_2^{}V_2^*&\\
&\sim_h V_1^{}\sigma\pi(\tilde{u})V_1^*+V_2^{}\sigma\pi(1_{\Td\ast_\C\Td})V_2^*&
\mathrm{by}\,\mathrm{assertion}\, 3)\\
&\sim_h \sigma\pi(1_{\Td\ast_\C\Td})
&\mathrm{by}\,\mathrm{assertion}\, 2)\\
&=1_{\Od\otimes C^*(\Z)}=\sigma(1_{\Od\ast_\C C^*(\Z)})
\end{array}$$ 
As a consequence, $\sigma(\ubf)\sim_h \sigma(s_1^{}s_1^*+s_2^{}\ubf s_2^*)$ 
in $\mathcal{U}(\Od\otimes C^*(\Z))$. 

\medskip\noindent5) The isomorphim $\Od\ast_{\C}C^*(\Z)\cong\Od\ast_\C\Od$ induces 
two embeddings $\sigma_1, \sigma_2$ of the \cst-algebra $\Od\ast_\C C^*(\Z)$ 
into the free product $\Od\ast_\C C^*(\Z)\ast_\C C^*(\Z)=\Od\ast_\C C^*(\mathbb{F}_2)$, 
where $\mathbb{F}_2$ is the free group with 2 generators. 
If the unitary $\mathring{u}$ is in $\mathcal{U}^0(\Od\ast_\C C^*(\Z))$, then the product 
$w:=\sigma_1(\mathring{u}) \sigma_2(\mathring{u}) \sigma_1(\mathring{u})^{-1}
 \sigma_2(\mathring{u})^{-1}$ belongs to $\mathcal{U}^0(\Od\ast_\C C^*(\mathbb{F}_2))$. 
We shall show in three steps that 
$w\not\in\mathcal{U}^0(\Od\ast_\C C^*(\mathbb{F}_2))$, 
and this will imply that the two unitaries $\ubf$ and $s_1^{}s_1^*+s_2^{}\ubf s_2^*$ 
are not homotopic in the compact group $\mathcal{U}(\Od\ast_\C C^*(\Z))$. 

Let $\beta$ be the action of the group $S^1=\mathbb{R}/\mathbb{Z}$ on the \cst-algebra $\Od$ 
given by $\beta_t(s_k)=e^{2\imath\pi t} s_k$ ($k=1, 2$). 
The subalgebra $A\subset\Od$ of $\beta$-invariant elements is the closure of 
the growing sequence of matrix \cst-algebras $A_n$ linearly generated by the elements 
${s_{\imath_1}^{}\ldots s_{\imath_n}^{}s_{\imath_{n+1}}^*\ldots s_{\imath_{2n}}^*}$ 
and where each $A_n$ embeds in $A_{n+1}$ by $a\mapsto s_1^{}as_1^*$. 
This sequence of monomorphisms extends to an endomorphism $\alpha$ on $A$. 
Call $\ddot{A}$ the inductive limit of the system 
$A\mathop{\rightarrow}\limits^{\alpha} A\mathop{\rightarrow}\limits^{\alpha}\ldots$ 
with corresponding embeddings $\mu_n: A\to\ddot{A}$ $(n\in\N$) and 
let $\ddot{\alpha}: \ddot{A}\to\ddot{A}$ be the unique automorphism 
satisfying $\ddot{\alpha}\left(\mu_n(a)\right)=\mu_n\left(\alpha(a)\right)$. 

\noindent\textit{Step 1. The unitary $w$ does not belong to the connected component 
$\mathcal{U}^0(A\ast_\C C^*(\Z))$.} \\
\noindent\textit{Proof}: Define for all integer $n\geq 2$ 
the compact space $U(n)$ of unitaries in $M_n(\C)$ and 
let $D_n$ be the \cst-algebra $D_n=C(U(n)\times U(n), M_n(\C))$. 
As noticed by B. Blackadar in Barcelona (\cite{Blac07}), 
a theorem by S. Araki, M. James and E. Thomas (\cite{AJT60}) 
implies that $w\not\in\mathcal{U}^0(A_n\ast_\C C^*(\Z ))$. 
Hence, $w\not\in\mathcal{U}^0(A\ast_\C C^*(\Z ))$ by passage to the limit 
since all the connecting maps $A_n\to A_{n+1}$ are injective. 

\smallskip
\noindent\textit{Step 2. The unitary $w$ in not in $\mathcal{U}^0(\ddot{A}\rtimes\Z \ast_\C C^*(\Z))$.} \\
\noindent\textit{Proof}: The density of the subset $\cup_n\mu_n(A)$ in $\ddot{A}$ implies that 
$w\not\in\mathcal{U}^0(\ddot{A}\ast_\C C^*(\Z))$. 
The epimorphism $f\mapsto f(1)$ from $C(S^1)=C^*(\Z)$ to $\C$ induces a surjection 
from $\ddot{A}\rtimes_{\ddot{\alpha}}\Z\ast_\C C^*(\Z)$ to the subalgebra $\ddot{A}\ast_\C C^*(\Z)$ 
and so $w$ could be in $\mathcal{U}^0(\ddot{A}\rtimes\Z \ast_\C C^*(\Z))$ only if 
$w$ was already in $\mathcal{U}^0(A\ast_\C C^*(\Z))$. 

\smallskip  
\noindent\textit{Step 3. The unitary $w$ does not belong to the connected component 
$\mathcal{U}^0(\Od\ast_\C C^*(\Z))$.} \\
\noindent\textit{Proof}: 
If $p\in\ddot{A}$ is the projection $p=\mu_0(1_A)$, 
then $p\left(\ddot{A}\rtimes_{\ddot{\alpha}}\Z\right)p\cong\Od$ 
(\cite[Subsection 2.1]{Cun77}, \cite[Claim 3.4]{DyShli01}). 
Hence the relation $w\not\in\mathcal{U}^0(\ddot{A}\rtimes\Z \ast_\C C^*(\Z))$ implies that 
$w\not\in\mathcal{U}^0(\Od\ast_\C C^*(\Z))$

\medskip\noindent6) Assertions 3) and 5) of the present Proposition imply that 
the $K_1$-trivial unitary $\tilde{u}$ is not in $\mathcal{U}^0(\Td\ast_\C\Td)$. 
Proposition~3.3 from \cite{Blan16} implies that 
the unital free product $\Td\ast_\C\Td$ is not $K_1$-injective. \qed 

\medskip 
\begin{rems}
1) 
The unital free product $\Td\ast_\C\Td$ does not have real rank $0$. 
Indeed, this would imply by Corollary 4.2.10 of \cite{Lin01} that 
the \cst-algebra $\Td\ast_\C\Td$ is K$_1$-injective. 

\noindent
2) The \cst-algebra $\Td$ is K$_1$-injective and so the amalgamated free product 
$\mathcal{U}(\Td)\mathop{\ast}\limits_{S^1}\mathcal{U}(\Td)$ embeds 
in the connected component $\mathcal{U}^0(\Td\ast_\C\Td)$\,. 
\end{rems} 
\medskip 

\begin{cor}\label{deformpropinf} 
There exists a unital continuous $C([0, 1])$-algebra with properly infinite fibres 
which is not a properly infinite \cst-algebra. 
\end{cor} 
\pf The \cst-algebra $\mathcal{D}\!:=\!\!\{f\in C([0, 1] , \Td\ast_\C\Td)\,;\, 
f(0)\in\jmath_0(\Td)\,\mathrm{and}\,f(1)\!\in\!\jmath_1(\Td)\,\}$ is 
a unital continuous $C([0, 1])$-algebra with properly infinite fibres. 
Proposition~3.3 of \cite{Blan16} and the above Proposition~\ref{prop2} 
imply that this \cst-algebra $\mathcal{D}$ is not properly infinite.

\href{mailto:Etienne.Blanchard@imj-prg.fr}{Etienne.Blanchard@imj-prg.fr} 

\address{IMJ-PRG, \quad UP7D - Campus des Grands Moulins, \quad Case 7012\\ 
${}$\quad F-75205 Paris Cedex 13} 
\end{document}